\RequirePackage{ifpdf}
\ifpdf 
\documentclass[pdftex]{sigma}
\else
\documentclass{sigma}
\fi

\begin{document}


\renewcommand{\thefootnote}{$\star$}

\renewcommand{\PaperNumber}{082}

\FirstPageHeading

\ShortArticleName{Higher Order Connections}

\ArticleName{Higher Order Connections\footnote{This paper is a
contribution to the Special Issue ``\'Elie Cartan and Dif\/ferential Geometry''. The
full collection is available at
\href{http://www.emis.de/journals/SIGMA/Cartan.html}{http://www.emis.de/journals/SIGMA/Cartan.html}}}

\Author{Michael G. EASTWOOD}

\AuthorNameForHeading{M.G. Eastwood}

\Address{Mathematical Sciences Institute,
Australian National University, ACT 0200, Australia}

\Email{\href{mailto:meastwoo@member.ams.org}{meastwoo@member.ams.org}}

\ArticleDates{Received March 12, 2009, in f\/inal form August 10, 2009;  Published online August 11, 2009}

\Abstract{The purpose of this article is to present the theory of higher order
connections on vector bundles from a viewpoint inspired by projective
dif\/ferential geometry.}

\Keywords{connections; jets; projective dif\/ferential geometry}

\Classification{53B05; 58A20}

\section{Introduction}\label{intro}
We begin with a few well-known remarks on commonplace linear connections. Let
$E$ denote a~smooth vector bundle on a~smooth manifold~$M$ (throughout this
article we work in the smooth category but all constructions go through {\it
mutatis mutandis\/} in the holomorphic category). A~{\it connection\/} on $E$
may be def\/ined as a splitting of the f\/irst jet exact sequence~\cite{s}
\[0\to\Lambda^1\otimes E\to J^1E\to E\to 0,\]
where $\Lambda^1$ is the bundle of $1$-forms on~$M$. Equivalently, a connection
on $E$ is a f\/irst order linear dif\/ferential operator
\[\nabla:E\to\Lambda^1\otimes E\]
whose symbol $\Lambda^1\otimes E\to\Lambda^1\otimes E$ is the identity. A
connection on $E$ induces a natural dif\/ferential operator
\[\nabla:\Lambda^1\otimes E\to\Lambda^2\otimes E,\enskip
\mbox{characterised by
$\nabla(\omega\otimes s)=d\omega\otimes s-\omega\wedge\nabla s$}\]
and the composition
\[E\xrightarrow{\,\nabla\,}\Lambda^1\otimes E\xrightarrow{\,\nabla\,}
\Lambda^2\otimes E\]
is a homomorphism of vector bundles called the {\it curvature\/} of~$\nabla$.

In this article, a {\it $k^{\mathit{th}}\!$ order connection\/} on $E$ is a
splitting of the $k^{\mathrm{th}}$ jet exact sequence
\begin{equation}\label{k-jet}
\textstyle 0\to\bigodot^k\!\Lambda^1\otimes E\to J^kE\to J^{k-1}E\to 0,
\end{equation}
where $\bigodot^k\!\Lambda^1$ is the $k^{\mathrm{th}}$ symmetric power
of~$\Lambda^1$. Libermann~\cite[p.~155]{l} considers these higher order
connections (more their {\it semi-holonomic\/} counterparts) but does not
pursue them so much along the lines done below. Other notions of higher order
connections are due to various authors including Ehresmann~\cite{e},
Virsik~\cite{v}, and Yuen~\cite{y}. A $k^{\mathrm{th}}$~order connection, as
above, is evidently equivalent to a $k^{\mathrm{th}}$ order linear dif\/ferential
operator
\[\textstyle\nabla^{(k)}:E\to \bigodot^k\!\Lambda^1\otimes E\]
whose symbol $\bigodot^k\!\Lambda^1\otimes E\to\bigodot^k\!\Lambda^1\otimes E$
is the identity. In the rest of this article, we motivate and extend the notion
of curvature {\it et cetera} to these higher order connections. On the way, we
shall encounter various useful constructions and assemble evidence for a
f\/inal conjecture.

\section{Interlude on projective geometry}\label{interlude}
Let $\partial_a$ denote the usual partial dif\/ferential operator
$\partial/\partial x^a$ on ${\mathbb{R}}^n$ with co\"ordinates~$x^a$. If a~smooth $1$-form $\omega_a$ is obtained as the exterior derivative $\partial_af$
of a smooth function~$f$, then $\partial_a\omega_b$ is necessarily symmetric in
its indices and, conversely, this condition is locally suf\/f\/icient to ensure
that $\omega_a=\partial_af$ for some~$f$. Otherwise said, if we denote the skew
part of a tensor by enclosing the relevant indices in square brackets, then
being in the kernel of the operator $\omega_a\mapsto\partial_{[a}\omega_{b]}$
is the local integrability condition for the range of $f\mapsto\partial_af$. Of
course, these operators are the f\/irst two in the de~Rham complex
\[\Lambda^0\xrightarrow{\,d\,}\Lambda^1\xrightarrow{\,d\,}\Lambda^2
\xrightarrow{\,d\,}\cdots\xrightarrow{\,d\,}\Lambda^{n-1}\xrightarrow{\,d\,}
\Lambda^n.\]

Now consider the dif\/ferential operator on~${\mathbb{R}}^n$
\[\textstyle d^{(k)}: \ \Lambda^0\to\bigodot^k\!\Lambda^1\enskip
\mbox{given by}\enskip
f\mapsto\partial_a\partial_b\cdots\partial_cf
\enskip\mbox{($k$ derivatives)}.\]
Evidently, if a symmetric covariant tensor $\omega_{ab\cdots c}$ with $k$
indices is of the form $\partial_a\partial_b\cdots\partial_cf$ for some $f$,
then $\partial_{[a}\omega_{b]c\cdots d}=0$. Just as in the case $k=1$, this
necessary condition is also locally suf\/f\/icient to identify the range
of~$d^{(k)}$. This is a result from projective dif\/ferential geometry, a~full
discussion of which may be found in~\cite{projectivenotes}. Here, suf\/f\/ice it
to give the following derivation. Let us
def\/ine a connection $\nabla_a$ on the bundle
${\mathbb{T}}\equiv\Lambda^0\oplus\Lambda^1$ by
\[\nabla_a\left[\begin{array}cf\\ \mu_b\end{array}\right]\equiv
\left[\begin{array}c\partial_af-\mu_a\\
\partial_a\mu_b\end{array}\right].\]
Notice that
\[\nabla_a\nabla_b\left[\begin{array}cf\\ \mu_c\end{array}\right]=
\nabla_a\left[\begin{array}c\partial_bf-\mu_b\\
\partial_b\mu_c\end{array}\right]=
\left[\begin{array}c\partial_a(\partial_bf-\mu_b)
-\partial_b\mu_a\\
\partial_a\partial_b\mu_c
\end{array}\right]=
\left[\begin{array}c\partial_a\partial_bf-\partial_a\mu_b
-\partial_b\mu_a\\
\partial_a\partial_b\mu_c\end{array}\right]\]
is symmetric in $ab$. In other words, the connection $\nabla_a$ on
${\mathbb{T}}$ is f\/lat. It follows immediately, that the coupled de~Rham
complex
\[\Lambda^0\otimes{\mathbb{T}}\xrightarrow{\,\nabla\,}
\Lambda^1\otimes{\mathbb{T}}
\xrightarrow{\,\nabla\,}\Lambda^2\otimes{\mathbb{T}}
\xrightarrow{\,\nabla\,}\cdots\xrightarrow{\,\nabla\,}
\Lambda^{n-1}\otimes{\mathbb{T}}
\xrightarrow{\,\nabla\,}\Lambda^n\otimes{\mathbb{T}}\]
is locally exact. In particular, we conclude that locally
\[\left[\begin{array}c\phi_a\\ \omega_{ab}\end{array}\right]=
\left[\begin{array}c\partial_af-\mu_a\\
\partial_a\mu_b\end{array}\right]\enskip\mbox{for some}\enskip
\left[\begin{array}cf\\ \mu_b\end{array}\right]
\enskip\mbox{if and ony if}\enskip
\left[\begin{array}c\partial_{[a}\phi_{b]}+\omega_{[ab]}\\
\partial_{[a}\omega_{b]c}\end{array}\right]=0.\]
In particular, if we take $\phi_a=0$ and $\omega_{ab}$ to be symmetric, then
this statement reads
\[\omega_{ab}=\partial_a\partial_bf\enskip\mbox{for some $f$ if and only if}
\enskip\partial_{[a}\omega_{b]c}=0,\]
as required in case $k=2$. The general case may be similarly derived from the
induced f\/lat connection on $\bigodot^{k-1}\!{\mathbb{T}}$. In order further to
untangle the consequences of the local exactness of these coupled de~Rham
sequences, it is useful to def\/ine various additional tensor bundles.
Let $\Theta^{p,q}$ be the bundle (used here on ${\mathbb{R}}^n$ but let us
maintain the same notation on a general
manifold) whose sections are
covariant tensors satisfying the following symmetries
\[\phi_{\mbox{\scriptsize$\underbrace{a\cdots b}_p$}}
{}_{\mbox{\scriptsize$\underbrace{cd\cdots e}_q$}}=
\phi_{[a\cdots b](cd\cdots e)}\enskip\mbox{such that}\enskip
\phi_{[a\cdots bc]d\cdots e}=0,\]
where enclosing indices in round brackets means to take the symmetric part.
These include the bundles we have encountered so far
\[\textstyle\Lambda^p=\Theta^{p,0}\quad\mbox{and}\quad
\bigodot^k\!\Lambda^1=\Theta^{1,k-1}\]
and also accommodate the local integrability conditions for the range of
$f\mapsto\partial_a\partial_b\cdots\partial_cf$. Sorting out the meaning of
local exactness for the coupled de~Rham complex
$\Lambda^\bullet\otimes\bigodot^{k-1}\!{\mathbb{T}}$, we f\/ind that there are
locally exact complexes on ${\mathbb{R}}^n$
\begin{equation}\label{specialBGG}\Lambda^0\xrightarrow{\,\partial^{(k)}\,}
\Theta^{1,k-1}\xrightarrow{\,\partial\,}\Theta^{2,k-1}
\xrightarrow{\,\partial\,}\cdots\xrightarrow{\,\partial\,}
\Theta^{n-1,k-1}\xrightarrow{\,\partial\,}
\Theta^{n,k-1}\end{equation}
for all $k\geq 1$ with the case $k=1$ being the de~Rham complex itself.
Details are left to the reader. These complexes are special cases of the
Bernstein--Gelfand--Gelfand (BGG) complex on real projective space
${\mathbb{RP}}_n$ viewed in a standard af\/f\/ine co\"ordinate patch
${\mathbb{R}}^n\hookrightarrow{\mathbb{RP}}_n$. Details may be found
in~\cite{projectivenotes}. An independent construction was given by
Olver~\cite{o}.

\section{Curvature}
The integrability conditions found in \S~\ref{interlude} provide the motivation
for the following construction.
\begin{theorem}\label{next} A $k^{\mathit{th}}\!$ order connection
$\nabla^{(k)}:E\to\bigodot^k\!\Lambda^1\otimes E=\Theta^{1,k-1}\otimes E$
canonically induces a~first order operator
$\nabla:\Theta^{1,k-1}\otimes E\to\Theta^{2,k-1}\otimes E$ characterised by
the following two properties
\begin{itemize}\itemsep=0pt
\item its symbol
$\Lambda^1\otimes\bigodot^k\!\Lambda^1\otimes E\to\Theta^{2,k-1}\otimes E$ is
$\delta\otimes{\mathrm{Id}}$ where
$\delta:\Lambda^1\otimes\bigodot^k\!\Lambda^1\to\Theta^{2,k-1}$ is the
tensorial homomorphism
$\phi_{abc\cdots d}\stackrel{\delta}{\longmapsto}\phi_{[ab]c\cdots d}$;
\item the composition
$E\xrightarrow{\,\nabla^{(k)}\,}\bigodot^k\!\Lambda^1\otimes E
\xrightarrow{\,\nabla\,}\Theta^{2,k-1}\otimes E$ has order $k-1$.
\end{itemize}
\end{theorem}
\begin{proof} If we choose an arbitrary local trivialisation of $E$ and local
co\"ordinates on $M$, then
\[s\stackrel{\nabla^{(k)}}
{\mapstochar\relbar\joinrel\relbar\joinrel\longrightarrow}
\overbrace{\partial_{(b}\partial_c\partial_d\cdots\partial_{e)}}^ks+
\Gamma_{bcd\cdots e}{}^{fg\cdots h}
\overbrace{\partial_f\partial_g\cdots\partial_h}^{k-1}s
+\mbox{lower order terms}\]
for a uniquely def\/ined tensor $\Gamma_{bcd\cdots e}{}^{fg\cdots h}$ symmetric
in both its lower and upper indices and having values in~${\mathrm{End}}(E)$.
But then
\[\omega_{bcd\cdots e}\stackrel{\nabla}{\longmapsto}
\partial_{[a}\omega_{b]cd\cdots e}
+\Gamma_{cd\cdots e[a}{}^{fg\cdots h}\omega_{b]fg\cdots h}\]
is forced by the two characterising properties of~$\nabla$.
\end{proof}
\begin{definition} We shall refer to the composition
\[E\xrightarrow{\,\nabla\circ\nabla^{(k)}\,}\Theta^{2,k-1}\otimes E\]
as the {\it curvature\/} of~$\nabla^{(k)}$. Of course, when $k=1$ this is the
usual notion of curvature $E\to\Lambda^2\otimes E$ for a (f\/irst order)
connection on~$E$. The operator $\partial^{(k)}$ on ${\mathbb{R}}^n$ has zero
curvature.
\end{definition}

An alternative construction, both of the operator
$\nabla:\bigodot^k\!\Lambda^1\otimes E\to\Theta^{2,k-1}\otimes E$ and the
curvature $E\to\Theta^{2,k-1}\otimes E$, may be
given by expressing higher order connections in terms of commonplace
connections on the jet bundle $J^{k-1}E$ as follows. Recall that the
{\it Spencer operator\/} is a canonically def\/ined f\/irst order dif\/ferential
operator ${\mathcal{S}}:J^\ell E\to\Lambda^1\otimes J^{\ell-1}E$ characterised
by the following properties~\cite[Propositions~4~and~5]{g}:
\begin{itemize}\itemsep=0pt
\item its symbol is
$\Lambda^1\otimes J^\ell E\xrightarrow{\,{\mathrm{Id}}\otimes\pi\,}
\Lambda^1\otimes J^{\ell-1}E$ where $\pi$ is the canonical jet projection;
\item the sequence
$E\xrightarrow{\,j^\ell\,}J^\ell E\xrightarrow{\,{\mathcal{S}}\,}
\Lambda^1\otimes J^{\ell-1}E$, where $j^\ell$ is the universal
$\ell^{\mathrm{th}}$ order dif\/ferential operator, is locally exact.
\end{itemize}
As a splitting of (\ref{k-jet}), we may regard a $k^{\mathrm{th}}$ order
connection as a homomorphism $h:J^{k-1}E\to J^k E$ such that $\pi\circ
h={\mathrm{Id}}$. Composing with the Spencer operator
\begin{equation}\label{create}
J^{k-1}E\xrightarrow{\,h\,}J^kE\xrightarrow{\,{\mathcal{S}}\,}
\Lambda^1\otimes J^{k-1}E\end{equation}
gives a f\/irst order dif\/ferential operator $\nabla\equiv{\mathcal{S}}\circ h$
whose symbol $\Lambda^1\otimes J^{k-1}E\to\Lambda^1\otimes J^{k-1}E$ is the
identity, in order words a connection on~$J^{k-1}E$. (In the holomorphic
category, Jahnke and Radlof\/f already observed~\cite{jr} that a splitting of the
jet exact sequence (\ref{k-jet}) implied the vanishing of the Atiyah
obstruction~\cite{a} to $J^{k-1}E$ admitting a connection but (\ref{create})
is stronger in actually creating the desired connection.)
\begin{theorem}\label{connections}
A $k^{\mathit{th}}\!$ order connection on a vector bundle $E$ induces a
commonplace connection on the jet bundle $J^{k-1}E$ with the following
properties
\begin{itemize}\itemsep=0pt
\item the composition
$J^{k-1}E\xrightarrow{\,\nabla\,}\Lambda^1\otimes J^{k-1}E
\xrightarrow{\,{\mathrm{Id}}\otimes\pi\,}\Lambda^1\otimes J^{k-2}E$ is the
Spencer operator;
\item its curvature $\kappa:J^{k-1}E\to\Lambda^2\otimes J^{k-1}E$ has values in
$\Lambda^2\otimes\bigodot^{k-1}\!\Lambda^1\otimes E
\hookrightarrow\Lambda^2\otimes J^{k-1}E$.
\end{itemize}
Conversely, a connection on $J^{k-1}E$ with these two properties uniquely
characterises a $k^{\mathit{th}}\!$ order connection on~$E$.
\end{theorem}
\begin{proof}
As observed in~\cite{g}, the Spencer operator induces f\/irst order dif\/ferential
operators
\begin{equation}\label{higherS}
{\mathcal{S}}:\Lambda^1\otimes J^\ell E\to\Lambda^2\otimes J^{\ell-1}E
\enskip\mbox{def\/ined by}\enskip
{\mathcal{S}}(\omega\otimes s)=
d\omega\otimes\pi s-\omega\wedge{\mathcal{S}}s\end{equation}
and there is a commutative diagram~\cite[(31)]{g}
\begin{equation}\label{diagram}\begin{array}{ccccccccc}
&&0&&E&=&E\\
&&\downarrow&&\makebox[0pt]{\scriptsize$j^k$}\downarrow\makebox[0pt]{}&&
\makebox[0pt]{\scriptsize$j^{k-1}$\quad}\downarrow\makebox[0pt]{}\\
0&\to&\bigodot^k\!\Lambda^1\otimes E&\to&J^kE&\xrightarrow{\,\pi\,}&
J^{k-1}E&\to&0\\
&&\downarrow&&
\makebox[0pt]{\scriptsize${\mathcal{S}}$}\downarrow\makebox[0pt]{}&&
\makebox[0pt]{\scriptsize${\mathcal{S}}$}\downarrow\makebox[0pt]{}\\
0&\to&\Lambda^1\otimes\bigodot^{k-1}\!\Lambda^1\otimes E&\to&
\Lambda^1\otimes J^{k-1}E&\xrightarrow{\,{\mathrm{Id}}\otimes\pi\,}&
\Lambda^1\otimes J^{k-2}E&\to&0\\
&&\downarrow&&
\makebox[0pt]{\scriptsize${\mathcal{S}}$}\downarrow\makebox[0pt]{}&&
\makebox[0pt]{\scriptsize${\mathcal{S}}$}\downarrow\makebox[0pt]{}\\
0&\to&\Lambda^2\otimes\bigodot^{k-2}\!\Lambda^1\otimes E&\to&
\Lambda^2\otimes J^{k-2}E&\xrightarrow{\,{\mathrm{Id}}\otimes\pi\,}&
\Lambda^2\otimes J^{k-3}E&\to&0\\
&&\vdots&&\vdots&&\vdots
\end{array}\end{equation}
with exact rows of vector bundle homomorphisms and locally exact columns of
linear dif\/ferential operators (apart from the f\/irst column, which consists of
homomorphisms starting with $-\iota\otimes{\mathrm{Id}}$ where
$\iota:\bigodot^k\!\Lambda^1\hookrightarrow
\Lambda^1\otimes\bigodot^{k-1}\!\Lambda^1$ is the natural inclusion).
For the f\/irst characterising property
of~$\nabla$, we compute from~(\ref{create}) and~(\ref{diagram}):
\[({\mathrm{Id}}\otimes\pi)\circ\nabla=
({\mathrm{Id}}\otimes\pi)\circ{\mathcal{S}}\circ h=
{\mathcal{S}}\circ\pi\circ h={\mathcal{S}}\circ{\mathrm{Id}}={\mathcal{S}},\]
as required. To compute the curvature $\kappa$ of $\nabla$ we must consider the
induced operator
\[\nabla:\Lambda^1\otimes J^{k-1}E\to\Lambda^2\otimes J^{k-1}E
\enskip\mbox{def\/ined by}\enskip
\nabla(\omega\otimes s)=d\omega\otimes s-\omega\wedge\nabla s.\]
Composing this formula with
${\mathrm{Id}}\otimes\pi:\Lambda^2\otimes J^{k-1}E\to\Lambda^2\otimes J^{k-2}E$
gives
\[({\mathrm{Id}}\otimes\pi)\circ\nabla(\omega\otimes s)=
d\omega\otimes\pi s-\omega\wedge({\mathrm{Id}}\otimes\pi)\circ\nabla s=
d\omega\otimes\pi s-\omega\wedge{\mathcal{S}}s,\]
by the f\/irst property of $\nabla$ established above. {From} (\ref{higherS}) we
conclude that
\begin{equation}\label{S}
({\mathrm{Id}}\otimes\pi)\circ\nabla={\mathcal{S}}:\Lambda^1\otimes J^{k-1}E
\to\Lambda^2\otimes J^{k-2}E.\end{equation}
Therefore,
\[({\mathrm{Id}}\otimes\pi)\circ\kappa=
({\mathrm{Id}}\otimes\pi)\circ\nabla\circ\nabla={\mathcal{S}}\circ\nabla=
{\mathcal{S}}\circ{\mathcal{S}}\circ h=0,\quad\mbox{because}
\enskip{\mathcal{S}}\circ{\mathcal{S}}=0.\]
{From} the exactness of
$0\to\Lambda^2\otimes\bigodot^{k-1}\!\Lambda^1\otimes E\to
\Lambda^2\otimes J^{k-1}E\to\Lambda^2\otimes J^{k-2}E\to 0$, it follows that
$\kappa$ takes values in $\Lambda^2\otimes\bigodot^{k-1}\!\Lambda^1\otimes E$,
as required.

Conversely, given a connection $\nabla$ on $J^{k-1}E$ satisfying the two
properties in the statement of the theorem, let us def\/ine a $k^{\mathrm{th}}$
order dif\/ferential operator
\[\nabla^{(k)}:E\to\Lambda^1\otimes J^{k-1}E\enskip\mbox{as the composition
$\nabla\circ j^{k-1}$}.\]
{From} the f\/irst property of $\nabla$ we observe that
\[({\mathrm{Id}}\otimes\pi)\circ\nabla^{(k)}=
({\mathrm{Id}}\otimes\pi)\circ\nabla\circ j^{k-1}={\mathcal{S}}\circ j^{k-1}
=0\]
and, with reference to~(\ref{diagram}), deduce that, in fact,
\[\textstyle
\nabla^{(k)}:E\to\Lambda^1\otimes\bigodot^{k-1}\!\Lambda^1\otimes E.\]
It is easy to check that its symbol $\bigodot^k\!\Lambda^1\otimes E\to
\Lambda^1\otimes\bigodot^{k-1}\!\Lambda^1\otimes E$ is the natural inclusion
$\iota\otimes{\mathrm{Id}}$. Therefore, to show that $\nabla^{(k)}$
is, in fact, a $k^{\mathrm{th}}$ order connection, it suf\/f\/ices to show that
$\nabla^{(k)}$ takes values in $\bigodot^k\!\Lambda^1\otimes E$ and,
with reference to~(\ref{diagram}), for this it suf\/f\/ices to show that
${\mathcal{S}}\circ\nabla^{(k)}=0$. For this we may compute using
(\ref{S}) and our def\/inition of~$\nabla^{(k)}$:
\[{\mathcal{S}}\circ\nabla^{(k)}=
({\mathrm{Id}}\otimes\pi)\circ\nabla\circ\nabla\circ j^{k-1}=
({\mathrm{Id}}\otimes\pi)\circ\kappa\circ j^{k-1}=
0\circ j^{k-1}=0,\]
as required. Finally, we must check that this construction does indeed provide
an inverse to setting $\nabla\equiv{\mathcal{S}}\circ h$. {From}
(\ref{diagram}), the usual splitting rigmarole produces
\[\textstyle
j^k-h\circ j^{k-1}:E\to\bigodot^k\!\Lambda^1\otimes E\hookrightarrow J^kE.\]
But, viewing via
$\bigodot^k\!\Lambda^1\otimes E
\stackrel{\iota\otimes{\mathrm{Id}}}
{\lhook\joinrel\relbar\joinrel\longrightarrow}
\Lambda^1\otimes\bigodot^{k-1}\!\Lambda^1\otimes E
\hookrightarrow\Lambda^1\otimes J^{k-1}E$ (as we were doing) gives
\[-{\mathcal{S}}\circ(j^k-h\circ j^{k-1})=
{\mathcal{S}}\circ h\circ j^{k-1}=\nabla\circ j^{k-1}.\]
Therefore, the combined ef\/fect of $\nabla^{(k)}\rightsquigarrow
h\rightsquigarrow\nabla\equiv{\mathcal{S}}\circ h\rightsquigarrow
\nabla^{(k)}\equiv\nabla\circ J^{(k-1)}$ is to end up back
where we started. To check that
$\nabla\rightsquigarrow\nabla^{(k)}\equiv\nabla\circ J^{(k-1)}
\rightsquigarrow h\rightsquigarrow\nabla\equiv{\mathcal{S}}\circ h$ is also
the identity is a similar unravelling of def\/initions and is left to the reader.
\end{proof}

Any construction starting with a $k^{\mathrm{th}}$ order connection
$\nabla^{(k)}:E\to\bigodot^k\!\Lambda^1\otimes E$ may, of course, be carried
out using a commonplace connection on the jet bundle $J^{k-1}E$ in accordance
with Theorem~\ref{connections}. Consider, for example, the composition
\begin{equation}\label{compo}\textstyle\bigodot^k\!\Lambda^1\otimes E
\stackrel{\iota\otimes{\mathrm{Id}}}
{\lhook\joinrel\relbar\joinrel\longrightarrow}
\Lambda^1\otimes\bigodot^{k-1}\!\Lambda^1\otimes E
\hookrightarrow\Lambda^1\otimes J^{k-1}E\xrightarrow{\,\nabla\,}
\Lambda^2\otimes J^{k-1}E\end{equation}
where
$\Lambda^1\otimes J^{k-1}E\xrightarrow{\,\nabla\,}\Lambda^2\otimes J^{k-1}E$
is the usual induced f\/irst order operator. By~(\ref{S}), if we further compose
with ${\mathrm{Id}}\otimes\pi:
\Lambda^2\otimes J^{k-1}E\to\Lambda^2\otimes J^{k-2}E$ then we obtain
\[\textstyle\bigodot^k\!\Lambda^1\otimes E
\stackrel{\iota\otimes{\mathrm{Id}}}
{\lhook\joinrel\relbar\joinrel\longrightarrow}
\Lambda^1\otimes\bigodot^{k-1}\!\Lambda^1\otimes E
\hookrightarrow\Lambda^1\otimes J^{k-1}E\xrightarrow{\,{\mathcal{S}}\,}
\Lambda^2\otimes J^{k-2}E,\]
which vanishes by dint of~(\ref{diagram}). Therefore (\ref{compo}) actually has
range in $\Lambda^2\otimes\bigodot^{k-1}\!\otimes E$ and, in fact, has range
in $\Theta^{2,k-1}\otimes E$ as follows. According to the commutative square
\[\begin{array}{ccc}
\Lambda^2\otimes\bigodot^{k-1}\!\Lambda^1\otimes E&\hookrightarrow&
\Lambda^2\otimes J^{k-1}E\\
\downarrow&&
\makebox[0pt]{\scriptsize${\mathcal{S}}$}\downarrow\makebox[0pt]{}\\
\Lambda^3\otimes\bigodot^{k-2}\!\Lambda^1\otimes E&\hookrightarrow&
\Lambda^3\otimes J^{k-2}E,
\end{array}\]
we must show that composing (\ref{compo}) with
${\mathcal{S}}:\Lambda^2\otimes J^{k-1}E\to\Lambda^3\otimes J^{k-2}E$ gives
zero. But
\[\begin{array}{ccc}
\Lambda^1\otimes J^{k-1}E&\xrightarrow{\,\nabla\,}&
\Lambda^2\otimes J^{k-1}E\\
\makebox[0pt]{$\kappa$}\downarrow\makebox[0pt]{}&&
\makebox[0pt]{\scriptsize${\mathcal{S}}$}\downarrow\makebox[0pt]{}\\
\Lambda^3\otimes J^{k-1}E&\xrightarrow{\,{\mathrm{Id}}\otimes\pi\,}&
\Lambda^3\otimes J^{k-2}E
\end{array}\]
also commutes and we see that ${\mathcal{S}}\circ\nabla=0$ from the curvature
restriction imposed by the second condition in Theorem~\ref{connections}. In
summary, the composition (\ref{compo}) takes values in
$\Theta^{2,k-1}\otimes E$ and, of course, it is the f\/irst order dif\/ferential
operator characterised in Theorem~\ref{next}. Similarly, the two conditions
imposed by Theorem~\ref{connections} on a connection on $J^{k-1}E$ imply that
its curvature
\[\kappa:J^{k-1}E\to\Lambda^2\otimes J^{k-1}E\]
actually takes values in $\Theta^{2,k-1}\otimes E$ and as a
$(k-1)^{\mathrm{st}}$ order dif\/ferential operator on $E$ it coincides with
$\nabla\circ\nabla^{(k)}$.

\section{Application to prolongation}
In joint work in progress with Rod Gover, higher order connections are used
as follows. Suppose $D:E\to F$ is a linear dif\/ferential operator of order $k$
with surjective symbol. The diagram
\begin{equation}\label{KT}\begin{array}{ccccccccc}
&&0&&0\\
&&\downarrow&&\downarrow\\
0&\to&K&\to&{\mathbb{T}}&\to&J^{k-1}E&\to&0\\
&&\downarrow&&\downarrow&&\|\\
0&\to&\bigodot^k\!\Lambda^1\otimes E&\to&J^kE&\to&J^{k-1}E&\to&0\\
&&\makebox[0pt]{\scriptsize$\sigma(D)\enskip\quad$}\downarrow\makebox[0pt]{}&&
\makebox[0pt]{\scriptsize$D$}\downarrow\makebox[0pt]{}\\
&&F&=&F\\
&&\downarrow&&\downarrow\\
&&0&&0
\end{array}\end{equation}
with exact rows and columns def\/ines the bundles $K$ and ${\mathbb{T}}$.
Furthermore, a splitting of
\begin{equation}\label{T}0\to K\to{\mathbb{T}}\to E\to 0\end{equation}
evidently splits the middle row of~(\ref{KT}). Thus, a splitting of (\ref{T})
induces a $k^{\mathrm{th}}$ order connection~$\nabla^{(k)}$ on~$E$ such that
$D=\sigma(D)\circ\nabla^{(k)}$. In combination with Theorem~\ref{connections},
we see that
\[D\phi=0\iff\nabla^{(k)}\phi=\omega\iff\nabla\tilde\phi=\omega,\]
where $\omega$ is some section of $K$ and
$\tilde\phi=j^{k-1}\phi\in\Gamma(M,J^{k-1}E)$. Now, according to
Theorem~\ref{next} and the discussion at the end of the previous section, the
operator $\nabla:\bigodot^k\!\Lambda^1\otimes E\to\Theta^{2,k-1}\otimes E$
applied to $\omega$ may be written as $\kappa\tilde\phi$ where $\kappa$ is the
curvature of the connection $\nabla$ on $J^{k-1}E$. The upshot of this
reasoning is that the equation $D\phi=0$ may be rewritten as the following
system
\begin{gather*}\nabla\tilde\phi = \omega,\\
\nabla\omega = \kappa\tilde\phi,
\end{gather*}
for $(\tilde\phi,\omega)$ a section of $J^{k-1}E\oplus K$. This is the f\/irst
step in prolonging the equation $D\phi=0$. With more care, this f\/irst step may
be taken more invariantly, ending up with a well-def\/ined f\/irst order
dif\/ferential operator on ${\mathbb{T}}$ independent of our choice of splitting
of~(\ref{T}). We shall see an example of this phenomenon in the following
section.

\section{Other BGG-like sequences}
In Theorem~\ref{next}, we saw the start of a sequence of dif\/ferential operators
modelled on the complex~(\ref{specialBGG}) from projective dif\/ferential
geometry. In fact, it is not too hard to extend this sequence all the way
\[E\xrightarrow{\,\nabla^{(k)}\,}\Theta^{1,k-1}\otimes E
\xrightarrow{\,\nabla\,}\Theta^{2,k-1}\otimes E
\xrightarrow{\,\nabla\,}\cdots\xrightarrow{\,\nabla\,}
\Theta^{n-1,k-1}\otimes E\xrightarrow{\,\nabla\,}
\Theta^{n,k-1}\otimes E\]
as a coupled version of~(\ref{specialBGG}). The ingredients for this
construction are the connection $\nabla$ on $J^{k-1}E$ from Theorem~\ref{next}
and its relation with Spencer operators ${\mathcal{S}}$ coming
from~(\ref{diagram}). Details are left to the reader.

Another BGG complex on ${\mathbb{R}}^n\hookrightarrow{\mathbb{RP}}_n$ starts
with the operator $\Lambda^1\to\bigodot^2\!\Lambda^1$ given by
\begin{equation}\label{flatKilling}\phi_a\longmapsto\partial_{(a}\phi_{b)}
\end{equation}
and continues with the second order operator (sometimes called the Saint Venant
operator)
\begin{equation}\label{SV}
h_{ab}\longmapsto\partial_a\partial_ch_{bd}
-\partial_b\partial_ch_{ad}
-\partial_a\partial_dh_{bc}
+\partial_b\partial_dh_{ac}.\end{equation}
This suggests that if we are given an arbitrary f\/irst order dif\/ferential
operator
\[D:\textstyle\Lambda^1\otimes E\to\bigodot^2\!\Lambda^1\otimes E\]
whose symbol is
\[\textstyle
\Lambda^1\otimes\Lambda^1\otimes E
\xrightarrow{\,\underbar\;\odot\underbar\;\otimes{\mathrm{Id}}\,}
\bigodot^2\!\Lambda^1\otimes E,\]
then there should be a canonically def\/ined second order operator
\[\textstyle \bigodot^2\!\Lambda^1\otimes E\to\Xi^{2,2}\otimes E\]
with the same symbol as~(\ref{SV}), where $\Xi^{p,q}$ is the bundle
whose sections are covariant tensors satisfying the following symmetries
\[\phi_{\mbox{\scriptsize$\underbrace{a\cdots b}_p$}}
{}_{\mbox{\scriptsize$\underbrace{cd\cdots e}_q$}}=
\phi_{[a\cdots b][cd\cdots e]}\enskip\mbox{such that}
\enskip \phi_{[a\cdots bc]d\cdots e}=0.\]
In fact, inspired by the full BGG-complex on ${\mathbb{RP}}_n$, we might expect
a coupled sequence
\begin{equation}\label{inspired}
\textstyle\Lambda^1\otimes E\xrightarrow{\,D\,}
\bigodot^2\!\Lambda^1\otimes E\to\Xi^{2,2}\otimes E\to
\Xi^{3,2}\otimes E\to\cdots\to\Xi^{n-1,2}\otimes E\to\Xi^{n,2}\otimes E.
\end{equation}
This is, indeed, the case. Although it is not clear what should be the
counterpart to Theorem~\ref{next}, we may canonically construct the desired
operators as follows. Consider what becomes of~(\ref{KT}):
\begin{equation}\label{specialKT}\begin{array}{ccccccccc}
&&0&&0\\
&&\downarrow&&\downarrow\\
0&\to&\Lambda^2\otimes E&\to&{\mathbb{T}}&\to&\Lambda^1\otimes E&\to&0\\
&&\downarrow&&\downarrow&&\|\\
0&\to&\Lambda^1\otimes\Lambda^1\otimes E&\to&J^1(\Lambda^1\otimes E)&\to&
\Lambda^1\otimes E&\to&0\\
&&\downarrow&&
\makebox[0pt]{\scriptsize$D$}\downarrow\makebox[0pt]{}\\
&&\bigodot^2\!\Lambda^1\otimes E&=&\bigodot^2\!\Lambda^1\otimes E\\
&&\downarrow&&\downarrow\\
&&0&&0
\end{array}\end{equation}
In particular, this diagram def\/ines ${\mathbb{T}}$ and also shows that a
splitting of
\[0\to\Lambda^2\otimes E\to{\mathbb{T}}\to\Lambda^1\otimes E\to 0\]
not only enables us to write sections of ${\mathbb{T}}$ as
\[\left[\begin{array}c\phi_a\\ \mu_{ab}\end{array}\right]\quad
\mbox{for}\enskip
\left\{\!\begin{array}{l} \phi_a\in\Gamma(\Lambda^1\otimes E),\\
\mu_{ab}\in\Gamma(\Lambda^2\otimes E),\end{array}\right.\]
but also splits the middle row of (\ref{specialKT}), i.e.\ def\/ines a
connection $\nabla_a$ on $\Lambda^1\otimes E$. In terms of this connection, the
operator $D$ is simply $\phi_a\stackrel{D}{\longmapsto}\nabla_{(a}\phi_{b)}$.
Now
consider the operator
\[{\mathbb{T}}\ni
\left[\begin{array}c\phi_a\\ \mu_{ab}\end{array}\right]\mapsto
\left[\begin{array}c\nabla_a\phi_b-\mu_{ab}\\
\nabla_{[a}\mu_{b]c}-\kappa_{abc}{}^d\phi_d
-\nabla_{[a}\mu_{c]b}+\kappa_{acb}{}^d\phi_d
-\nabla_{[b}\mu_{c]a}+\kappa_{bca}{}^d\phi_d
\end{array}\right]\in\Lambda^1\otimes{\mathbb{T}},\]
where $\kappa:\Lambda^1\otimes E\to\Lambda^2\otimes\Lambda^1\otimes E$ is the
curvature of~$\nabla$. It is a connection on ${\mathbb{T}}$ and a tedious
computation verif\/ies that it is independent of choice of splitting
of~${\mathbb{T}}$. Now consider the coupled de~Rham sequence with values
in~${\mathbb{T}}$ derived from this connection
\begin{equation}\label{coupled}\begin{array}{ccccccccc}
{\mathbb{T}}&\xrightarrow{\,\nabla\,}&\Lambda^1\otimes{\mathbb{T}}
&\xrightarrow{\,\nabla\,}&\Lambda^2\otimes{\mathbb{T}}
&\xrightarrow{\,\nabla\,}&\Lambda^3\otimes{\mathbb{T}}
&\xrightarrow{\,\nabla\,}&\cdots\\
\|&&\|&&\|&&\|\\
\Lambda^1\otimes E
&&\!\!\Lambda^1\otimes\Lambda^1\otimes E\!\!
&&\!\!\Lambda^2\otimes\Lambda^1\otimes E\!\!
&&\!\!\Lambda^3\otimes\Lambda^1\otimes E\!\!\\
\oplus&\nearrow&\oplus&\nearrow&\oplus&\nearrow&\oplus\\
\Lambda^2\otimes E
&&\!\!\Lambda^1\otimes\Lambda^2\otimes E\!\!
&&\!\!\Lambda^2\otimes\Lambda^2\otimes E\!\!
&&\!\!\Lambda^3\otimes\Lambda^2\otimes E,\!\!
\end{array}\end{equation}
noticing that the restricted operators $\nearrow$ are simply homomorphisms
given by
\[\mu_{\mbox{\scriptsize$\underbrace{a\cdots b}_p$}}{}_{cd}\mapsto
-\mu_{[a\cdots bc]d}.\]
When $p=0$ this homomorphism is injective. When $p=1$ it is an isomorphism. For
$p\geq 2$ it is surjective with $\Xi^{p,2}\otimes E$ as kernel. It is now just
diagram chasing to extract (\ref{inspired}) from (\ref{coupled}).

The Killing operator in Riemannian geometry provides a good example of a f\/irst
order linear dif\/ferential operator to which the reasoning above may be
applied. In this example, the bundle~$E$ is trivial and
\[\textstyle D:\Lambda^1\to\bigodot^2\!\Lambda^1\quad\mbox{is given by }
\phi_a\mapsto\nabla_{(a}\phi_{b)},\]
where $\nabla_a$ is the Levi-Civita connection. It is a straightforward
generalisation of (\ref{flatKilling}). Similarly, the f\/lat operator (\ref{SV})
is modif\/ied by replacing $\partial_a$ by $\nabla_a$ but also by adding suitable
zeroth order curvature terms. Details may be found in~\cite{projectivenotes}.
Both of these dif\/ferential operators have geometric interpretations. The
Killing operator itself gives the inf\/initesimal change in the Riemannian metric
$g_{ab}$ due to the f\/low of a vector f\/ield~$\phi^a$. The next operator
\[\textstyle\bigodot^2\!\Lambda\xrightarrow{\,\nabla^{(2)}\,}\Xi^{2,2}\]
gives the inf\/initesimal change in the Riemann curvature tensor due to a
perturbation of $g_{ab}$ by an arbitrary symmetric tensor (i.e.\ replace
$g_{ab}$ by $g_{ab}+\epsilon h_{ab}$ for suf\/f\/iciently small~$\epsilon$,
compute the Riemannian curvature for this new metric, dif\/ferentiate
in~$\epsilon$, and then set $\epsilon=0$). The next operator is an
inf\/initesimal manifestation of the Bianchi identity. This particular BGG
complex on ${\mathbb{RP}}_n$
\[\textstyle\Lambda^1\xrightarrow{\,\nabla\,}\bigodot^2\!\Lambda^1
\xrightarrow{\,\nabla^{(2)}\,}\Xi^{2,2}\xrightarrow{\,\nabla\,}
\Xi^{3,2}\xrightarrow{\,\nabla\,}\cdots\xrightarrow{\,\nabla\,}
\Xi^{n-1,2}\xrightarrow{\,\nabla\,}\Xi^{n,2}\]
was also constructed by Calabi~\cite{calabi} as the Riemannian deformation
complex for the constant curvature metric (only constant curvature metrics are
projectively f\/lat).

In three dimensions, the deformation of a Riemannian metric coincides with
the mathematical formulation of elasticity in continuum mechanics (see,
e.g.~\cite{c}). In three dimensions, we may also choose a volume form
$\epsilon_{abc}$ to ef\/fect an isomorphism
$\Xi^{2,2}\cong\bigodot^2\!\Lambda^1$ (a ref\/lection of the fact that in three
dimensions there is only Ricci curvature) and rewrite (\ref{SV}) as
\[\textstyle\bigodot^2\!\Lambda^1\ni h_{ab}\longmapsto
\epsilon_a{}^{cd}\epsilon_b{}^{ef}\partial_c\partial_eh_{df}\in
\bigodot^2\!\Lambda^1\quad
\mbox{(sometimes written as $h\mapsto\mathrm{curl\,curl}\,h$)}.\]
Also $\Xi^{3,2}\cong\Lambda^1$ and the {\it linearised elasticity complex\/}
becomes
\[\textstyle\Lambda^1\xrightarrow{\,\nabla\,}\bigodot^2\!\Lambda^1
\xrightarrow{\,\nabla^{(2)}\,}\bigodot^2\!\Lambda^1\xrightarrow{\,\nabla\,}
\Lambda^1,\]
usually interpreted as
${\mathit{displacement}}\mapsto{\mathit{strain}}\mapsto{\mathit{stress}}
\mapsto{\mathit{load}}$. A derivation of the complex in this form
on~${\mathbb{RP}}_3$ (by means of a coupled de~Rham complex as above) is given
in~\cite{srni}. The close link between BGG complexes and coupled de~Rham
complexes (in the f\/lat case) has recently been modif\/ied and then used by
Arnold, Falk, and Winther~\cite{afw} to give new stable f\/inite element schemes
applicable to numerical elasticity.

Having seen two examples thereof, it is natural to conjecture that there are
canonically def\/ined sequences of dif\/ferential operators modelled on the
general projective BGG complex. However, this remains a conjecture. Notice
that there is no direct link between projective dif\/ferential geometry and the
constructions in this article (and we are not using that the Killing operator
considered above happens to be projectively invariant when suitably
interpreted~\cite{projectivenotes}). More challenging cases of this conjecture
are to ask, for $s\geq 2$, if
\[D:\textstyle\bigodot^s\!\Lambda^1\otimes E\to
\bigodot^{s+1}\!\Lambda^1\otimes E\]
is an arbitrary f\/irst order dif\/ferential operator whose symbol is
\[\textstyle
\Lambda^1\otimes\bigodot^s\!\Lambda^1\otimes E
\xrightarrow{\,\underbar\;\odot\underbar\;\otimes{\mathrm{Id}}\,}
\bigodot^{s+1}\!\Lambda^1\otimes E,\]
whether there is a canonically def\/ined $(s+1)^{\mathrm{st}}$ order operator
\[\textstyle \bigodot^{s+1}\!\Lambda^1\otimes E\to
\Xi^{s+1,s+1}\otimes E\]
whose symbol is ${\mathrm{proj}}\otimes{\mathrm{Id}}$, where
\[\textstyle \bigodot^{s+1}\!\Lambda^1\otimes\bigodot^{s+1}\!\Lambda^1\to
\Xi^{s+1,s+1}\]
is induced by the canonical projection of
${\mathrm{GL}}(n,{\mathbb{R}})$-modules
\[\raisebox{5pt}{$\begin{picture}(48,8)
\put(0,0){\line(1,0){48}}
\put(0,8){\line(1,0){48}}
\put(0,0){\line(0,1){8}}
\put(8,0){\line(0,1){8}}
\put(16,0){\line(0,1){8}}
\put(28,4){\makebox(0,0){$\cdots$}}
\put(40,0){\line(0,1){8}}
\put(48,0){\line(0,1){8}}
\end{picture}\otimes
\begin{picture}(48,8)
\put(0,0){\line(1,0){48}}
\put(0,8){\line(1,0){48}}
\put(0,0){\line(0,1){8}}
\put(8,0){\line(0,1){8}}
\put(16,0){\line(0,1){8}}
\put(28,4){\makebox(0,0){$\cdots$}}
\put(40,0){\line(0,1){8}}
\put(48,0){\line(0,1){8}}
\end{picture}\to\;$}
\begin{picture}(48,16)
\put(0,0){\line(1,0){48}}
\put(0,8){\line(1,0){48}}
\put(0,16){\line(1,0){48}}
\put(0,0){\line(0,1){16}}
\put(8,0){\line(0,1){16}}
\put(16,0){\line(0,1){16}}
\put(28,4){\makebox(0,0){$\cdots$}}
\put(28,12){\makebox(0,0){$\cdots$}}
\put(40,0){\line(0,1){16}}
\put(48,0){\line(0,1){16}}
\end{picture}\;.\]
Even restricting attention to f\/irst order operators~$D$, there are many more
examples of this conjecture that might be considered. In general, the
conjecture applies to operators whose symbol is induced by a {\it Cartan
product} as in~\cite{bceg}.

\subsection*{Acknowledgements}
It is a pleasure to acknowledge many very useful discussions with Rod Gover.
Support from the Australian Research Council is also gratefully acknowledged.

\pdfbookmark[1]{References}{ref}
\LastPageEnding


\begin{thebibliography}{99}

\footnotesize\itemsep=0pt

\bibitem{afw}
Arnold D.N., Falk R.S.,  Winther R.,
Finite element exterior calculus, homological techniques, and applications,
{\it Acta Numer.} {\bf 15} (2006), 1--155.

\bibitem{a}
Atiyah M.F.,
Complex analytic connections in f\/ibre bundles,
{\it Trans. Amer. Math. Soc.} {\bf 85} (1957), 181--207.

\bibitem{bceg}
Branson T.P., \v{C}ap A., Eastwood M.G., Gover A.R.,
Prolongations of geometric overdetermined systems,
{\it Internat. J. Math.} {\bf 17} (2006), 641--664,
\href{http://arxiv.org/abs/math.DG/0402100}{math.DG/0402100}.


\bibitem{calabi}
Calabi E.,
On compact Riemannian manifolds with constant curvature I,
in Dif\/ferential Geometry,
Proc. Sympos. Pure Math., Vol.~III,
Amer. Math. Soc., Providence, R.I.,  1961, 155--180.

\bibitem{c}
Ciarlet P.G.,
An introduction to dif\/ferential geometry with application to elasticity,
{\it J. Elasticity} {\bf 78/79} (2005), no. 1-3, 1--215.

\bibitem{srni}
Eastwood M.G.,
A complex from linear elasticity,
in The Proceedings of the 19th Winter School ``Geometry and Physics'' (Srni, 1999),
{\it Rend. Circ. Mat. Palermo (2) Suppl.}  (2000), no.~63, 23--29.


\bibitem{projectivenotes}
Eastwood M.G.,
Notes on projective dif\/ferential geometry,
in Symmetries and Overdetermined Systems of Partial Dif\/ferential Equations,
{\it IMA Vol. Math. Appl.}, Vol.~144, Springer, New York, 2008,  41--60.

\bibitem{e}
Ehresmann  C.,
Connexions d'ordre sup\'erieur,
in Atti Quinto Congr. Un. Mat. Ital. (Pavia-Torino, 1955),
Edizioni Cremonese, 1956,  326--328.

\bibitem{g}
Goldschmidt H.,
Prolongations of linear partial dif\/ferential equations. I.~A conjecture of \'Elie Cartan,
{\it Ann. Sci. \'Ecole Norm. Sup. (4)} {\bf 1} (1968), 417--444.

\bibitem{jr}
 Jahnke P., Radlof\/f I.,
Splitting jet sequences,
{\it Math. Res. Lett.} {\bf 11} (2004), 345--354,
\href{http://arxiv.org/abs/math.AG/0210454}{math.AG/0210454}.

\bibitem{l}
Libermann P.,
Sur la g\'eom\'etrie des prolongements des espaces f\/ibr\'es vectoriels,
{\it  Ann. Inst. Fourier (Grenoble)} {\bf 14} (1964), 145--172.

\bibitem{o}
Olver P.J.,
Dif\/ferential hyperforms I,
Mathematics Report 82-101, University of Minnesota, 1982,
available at \url{http://www.math.umn.edu/~olver/a_/hyper.pdf}.

\bibitem{s}
Spencer D.C.,
Overdetermined systems of linear partial dif\/ferential equations,
{\it Bull. Amer. Math. Soc.} {\bf 75} (1969), 179--239.

\bibitem{v}
Virsik J.,
Non-holonomic connections on vector bundles.~II,
{\it Czechoslovak Math. J.} {\bf 17} (1967), 200--224.

\bibitem{y}
Yuen P.C.,
Higher order frames and linear connections,
{\it Cahiers Topologie G\'eom. Diff\'erentielle} {\bf 12} (1971), 333--371.

\end{thebibliography}
\end{document}